\documentclass[11pt, reqno]{amsproc}
\usepackage{amssymb}
\usepackage{amsmath}
\usepackage{mathdots}
\usepackage{amsbsy}
\usepackage{amscd}
\usepackage{amsthm}
\usepackage{xcolor}

\usepackage{url}
\usepackage[colorlinks,citecolor=red,urlcolor=blue,
bookmarks=false,hypertexnames=true]{hyperref}

\begin{document}
\title[
Comment on: The groups of order $p^6$ ($p$ an odd prime)]
{Comment on: \\ The groups of order $p^6$ ($p$ an odd prime) \\
By Rodney James, 
{\it Math.\ Comput.\ }{\bf 34} (1980), 613--637}

\begin{abstract}
Over the years various errors have been identified in the 1980 paper of James on the
groups of order $p^6$, where $p$ is an odd prime. 
Here we summarise them. 
\end{abstract}

\author{E.A.\ O'Brien}
\address{Department of Mathematics, University of Auckland, Auckland, New Zealand}
\email{e.obrien@auckland.ac.nz}
\author{Sunil Kumar Prajapati}
\address{Indian Institute of Technology Bhubaneswar, Arugul Campus, Jatni,
Khurda-752050, India}
\email{skprajapati@iitbbs.ac.in}
\author{Ayush Udeep}
\address{Indian Institute of Technology Bhubaneswar, Arugul Campus,
        Jatni, Khurda-752050, India}
\email{udeepayush@gmail.com}

\maketitle

\footnote{
We thank M.F. Newman for his assistance in preparing this report.
This work was supported in part by the Marsden Fund
of New Zealand via grant 20-UOA-2030.

}

\section{Introduction}

Over the years errors have been identified in the 
1980 paper \cite{James80} of James on the
groups of order $p^6$. 

There is an independent determination of these groups \cite{p6}.
It differs significantly from that in \cite{James80}.
For the fully regular case, $p \ge 7$, there are
\begin{center}
$3p^2+39p+344+24 \gcd(p-1,3)+11 \gcd(p-1,4)+2\gcd(p-1,5)$
\end{center}
isomorphism types.
The formula also applies for $p=5$.
For $p=3$ there are 504 isomorphism types.
Presentations of the groups for $p\geq 7$, 
organised by isoclinism family,
are explicitly recorded in \cite{Arxiv}. 
An electronic version which generates the presentations 
in {\sc Magma} \cite{Magma} is available at \cite{github}.
For all odd primes the groups are available
electronically as part of {\sc SmallGroups} \cite{Small} distributed with 
both {\sf GAP} \cite{GAP4} and {\sc Magma}. 
We record here some of the more significant errors we are aware of.
We make no claim that our list is comprehensive.

Girnat \cite{Girnat} lists presentations for the groups of
order dividing $p^5$ for primes $p > 3$; 
presentations of these (and the corresponding 3-groups) 
are available electronically as part of {\sc SmallGroups}. 

\section{Some errors}
\subsection{ $p \ge 5$}
We list the errors according to isoclinism family.

\begin{itemize}

\item Family $\Phi_{2}$\\
The group
$\Phi_{2}(1^6) = \Phi_{2}(1^5) \times (1)$.

\item Family $\Phi_{8}$\\
The presentations for the groups $\Phi_{8}(42)$ and 
$\Phi_{8}(33)$ are the same.

\item Family $\Phi_{12}$\\
The number of groups in this family is $p+13$ and not 15. 

\item Family $\Phi_{15}$\\
The number of groups in this family is $p+3$ and not 
$\frac {1}{4} \, \left ( p^2 + p + 10 + (p-1, 4) \right )$. 

\item Family $\Phi_{21}$\\
There are many problems: both missing isomorphism types
and duplicate copies. 
For example, there are 90 groups of order $7^6$ in this family;
the list contains 92 presentations and 68 isomorphism types.
There are 282 groups of order $13^6$; the list contains
286 presentations and 211 isomorphism types.

\item Family $\Phi_{26}$\\
The presentations given for ${\Phi_{25+x}}(222){a_r}$ in the note ``Added
in Proof" based on the work of K\"{u}pper \cite{Kupper} are valid 
only for $\Phi_{25}$.  

\item Family $\Phi_{30}$\\
For each prime one group is missing from this family.

\item Families $\Phi_{31}$ and $\Phi_{32}$ \\
The  definitions of $y$ and $j$ are inaccurate.
Family $\Phi_{31}$ has seven groups and $\Phi_{32}$ five groups. 

\item Family $\Phi_{43}$\\
The generator $\gamma$ should be added to the list of generators.  
The defining equation listed for $k$ and $\ell$ is incorrect. 

\end{itemize}

\subsection{ $p = 3$}
James reports that there are 504 groups, but 
his list has many errors. 
These include a group missing from Family 28; 
isomorphisms within Family 10 and Family 43; and 
isomorphisms between groups in Family 40 and Family 41. 
Some of the information in Table 4.1 does not apply:
for example, the derived groups 
in Family 25 and Family 26 have type (111). 

\end{document}